\documentclass[11pt]{amsart}
\usepackage{amsmath, latexsym, amssymb, times}
\usepackage[colorlinks=true,pagebackref,hyperindex]{hyperref}
\newtheorem{theorem}{Theorem}[section]
\newtheorem{lemma}[theorem]{Lemma}
\newtheorem{prop}[theorem]{Proposition}
\newtheorem{cor}[theorem]{Corollary}

\newtheorem{notn}[theorem]{Notation and Discussion}
\theoremstyle{definition}
\newtheorem{example}[theorem]{Example}

\newcommand{\R}{\mathbb{R}}
\newcommand{\Z}{\mathbb{Z}}
\newcommand{\N}{\mathbb{N}}

\newcommand{\e}{\mathbf{e}}
\newcommand{\lam}{\ensuremath{{\lambda}}}

\newcommand{\hht}{\mathrm{ht}}

\newcommand{\codim}{\mathrm{codim}}

\newcommand{\xvec}[1]{\ensuremath{x_{1},\ldots,x_{#1}}}
\newcommand{\lcm}{\ensuremath{\mathrm{lcm}}}
\newcommand{\grp}{\ensuremath{\mathrm{grp}}}
\newcommand{\pos}{\ensuremath{\mathrm{pos}}}

\newcommand{\blam}{\boldsymbol{\lambda}}
\newcommand{\bomega}{\boldsymbol{\omega}}
\newcommand{\bbeta}{\boldsymbol{\beta}}

\newcommand{\balpha}{\boldsymbol{\alpha}}
\newcommand{\bgamma}{\boldsymbol{\gamma}}

\newcommand{\bmu}{\boldsymbol{\mu}}

\newcommand{\omvec}[1]{\ensuremath{\omega_{1},\ldots,\omega_{#1}}}
\newcommand{\M}{\mathfrak{m}}
\newcommand{\mcr}{\mathcal{R}}

\begin{document}
\title{Serre's Condition R$_{\ell}$ for Affine Semigroup Rings}
\author{Marie A. Vitulli}
\address{Department of Mathematics, 1222-University of Oregon, Eugene, OR  97403-1222}
\email{vitulli@uoregon.edu}
\thanks{}

 

 \begin{abstract} In this note we characterize the affine semigroup rings $K[S]$ over an arbitrary field  $K$ that satisfy condition $\mathrm{R}_{\ell}$ of Serre.  Our characterization is in terms of the face lattice of the positive cone $\pos(S)$ of $S$. We start by reviewing some basic facts about the faces of $\pos(S)$ and consequences for the monomial primes of $K[S]$.   After proving our characterization we turn our attention to the Rees algebras of a special class of monomial ideals in a polynomial ring over a field.  In this special case, some of the characterizing criteria are always satisfied.  We give examples of nonnormal affine semigroup rings that satisfy $\mathrm{R}_2$.  
 \end{abstract}
 
 \maketitle
 
\section{Introduction}

\medskip
The class of affine semigroup rings is rich with examples that combine the flavors of convex geometry and commutative algebra.  The structure of the semigroup ring $K[S]$ is
 intimately related to the structure of the affine semigroup $S$ and the cone $\pos(S)$ spanned 
 by $S$.  For example, it is well known that $K[S]$ is normal if and only if $S$ contains all integral points
  of $\pos(S)$ (see \cite{BH}).   Normal affine semigroup rings are Cohen-Macaulay by a theorem of 
  Hochster \cite{Hoch}.  Ishida \cite{Ish} characterized the S$_2$-ness of $K[S]$ in terms of $S$ and
   the facets of $\pos(S)$.   In \cite{GW1} Goto and Watanabe announced a characterization of those affine semigroups $S$ for which  $K[S]$ is Cohen-Macaulay; Hoa and Trung gave a corrected characterization  in terms of both $S$ and  the cone spanned by $S$ over the rational numbers in \cite{HT} and \cite{HT2}.  
  In an earlier paper \cite{Vi} the author characterized those affine semigroup rings 
   which satisfy Serre's condition R$_1$.   In this note we characterize those affine semigroup
    rings $K[S]$ over an arbitrary field  $K$ which satisfy condition $\mathrm{R}_{\ell}$
     of Serre.  Our characterization is in terms of the face lattice of the positive cone
      $\pos(S)$ of $S$. We start by recalling some basic facts about the faces of
       $\pos(S)$ and consequences for the monomial primes of $K[S]$.   After proving our
        characterization we turn our attention to the Rees algebras of a special class of
         monomial ideals in a polynomial ring over a  field.  We may view these as affine semigroup rings;  the associated affine semigroups are in some sense complementary to the class of polytopal semigroups introduced and studied by Bruns,  Gubeladze, and Trung in \cite{BG}  and \cite{BGT}. In this special case, most
          of the characterizing criteria are always satisfied.  We give examples of nonnormal
           affine semigroup rings that satisfy $\mathrm{R}_2$.  
           
           For the fundamentals on convex geometry we refer the reader to \cite{Ew}, \cite{Grun}, or \cite{Zieg}. For background on monoids and semigroup rings one can consult \cite{Gil}.
We make the standard assumptions that an affine semigroup $S$ is a subsemigroup of $\Z^n$ and that 
$\grp(S) = \Z^n$ for some positive integer $n$.   Consider the \textit{positive cone} $\pos(S)$ of $S$ defined by
$$\pos(S) = \{ c_1\balpha_1 + c_2\balpha_2 + \cdots + c_m\balpha_m \mid m \ge 0, \balpha_i \in S,  c_i \in \R_{\ge} \; (i = 1, \ldots, m)\},$$
where $\R_{\ge}$ denotes the set of nonnegative real numbers. Recall that $\pos(S)$ is a  polyhedral cone, that is, $\pos(S)$ is the intersection of finitely many positive halfspaces    $H_i^+ = \{ \balpha \in \R^n \mid \sigma_i(\balpha) \ge 0\}$, where $\sigma_i$ is a 
   linear form on $\R^n$.   Since $S$ is a finitely generated subsemigroup of $\Z^n$ we may 
   assume that each $\sigma_i$ has rational coefficients, that is, $\pos(S)$ is a \textit{rational polyhedral
    cone}.  After scaling we may assume that  the coefficients of each $\sigma_i$ are relatively prime
     integers; we call such a \textit{primitive linear form}.  Recall that a \textit{supporting hyperplane} of $\pos(S)$ is a hyperplane $H$ such that  $\pos(S) \cap H \ne \emptyset$ and $\pos(S)$ lies on one
       side of $H$.  A \textit{face} of $\pos(S)$ is the intersection of $\pos(S)$ and a supporting hyperplane.
         The faces of $\pos(S)$ are again rational polyhedral cones. By the \textit{dimension} $\dim(F)$ of $F$ we mean the dimension of the vector space spanned by $F$ and by the \textit{codimension} $\codim(F)$ we mean $n - \dim(F)$.  A face of codimension one is called a \textit{facet}. If we represent  a polyhedral cone $C$ as an irredundant intersection of positive halfspaces  $C= \cap_{i = 1}^r H_i^+$  then $F_1 = C \cap H_1, \ldots, F_r = C \cap H_r$ are the facets of $C$ (see \cite[Section 2.6]{Grun} ).
  
By a \textit{monomial }in $\mcr=K[S]$ we mean an element of the form $x^{\balpha}$ and by a \textit{monomial ideal  }we mean an ideal generated by monomials. There is an order-reversing bijective correspondence between the nonempty faces of $\pos(S)$ and the monomial primes of $\mcr$.  Indeed, the monomial prime of $\mcr$ corresponding to the nonempty face $F$ of $\pos(S)$ is  $P_F = (x^{\balpha} \mid \balpha \in S \setminus F )$ (e.g. see \cite{BH}).  We let $\M$ denote the ideal of $K[S]$ generated by all noninvertible monomials; it is the maximal monomial prime of $K[S]$.   Notice that we are considering the zero ideal to be monomial.  Finally, we let $\e_1, \ldots , \e_n$ denote the standard basis vectors for $\R^n$. 

The maximal proper faces of a polyhedral cone $C$ are precisely the facets of $C$. If $P = P_F$ is the monomial prime of height $d$ in an affine semigroup ring $K[S]$ corresponding to the face $F$ of $\pos(S)$, then there exists a chain of monomials primes of length $d$ descending from $P$ (see \cite[Prop. 1.2.1]{GW2} ) and $\hht(P) = \codim(F)$.

\section{Serre's Regularity Condition for Affine Semigroup Rings} 

We start this section with a basic result about $\Z^n$-graded rings, which is crucial for our purposes. A version for $\Z$-graded rings is well known (e.g. see  \cite[Ex. 2.24]{BH}).  Then we specialize to the case of an affine semigroup ring defined over a field.  Recall that if $P$ is a prime ideal of a $Z^n$-graded ring $R$ then $P^*$ denotes the largest homogeneous ideal of $R$ that is contained in $P$ and $R_{(P)}$  the homogeneous localization of $R$ at $P$, i.e., $R_{(P)} = S^{-1}R$, where $S$ is the set of homogeneous elements of $R$ that are not in $P$.  The graded ring $R_{(P)}$ is an example of a ${}^*$local ring, that is, a graded ring with a unique maximal homogeneous ideal.

\begin{prop}  \label{prop regular graded ring} Let $R$ be a Noetherian $\Z^n$-graded ring.   Then, $R$ is  regular  if and only if $R_{(P)}$ is regular, for every homogeneous prime ideal $P$ of $R$.  Furthermore, for a ${}^*$local ring $R$ with unique maximal homogeneous ideal $\M$, the ring $R$ is regular if and only if $R_{\M}$ is regular.
\end{prop}
\begin{proof}   First suppose that $R$ is regular.  Let $P$ be a homogeneous prime ideal.  Then, $R_P$ is a regular local ring.  We must show that the ${}^*$local ring $\mcr := R_{(P)}$ is regular.  If $P\mcr$ is a maximal 
ideal of $\mcr$ then $R_{(P)} = R_P$ is a regular local ring.  So assume $P\mcr$ is not maximal and
 choose a prime ideal $Q$ of $R$ such that $Q^*= P$.  By assumption, $R_P = R_{Q^*} \cong \mcr_{Q^*\mcr}$ is a regular local ring. 
  Let  $\mathcal{P}$ be any prime ideal of $\mcr$.  Then 
  $\mathcal{P}^* \subseteq Q^*\mcr  \Rightarrow \mcr _{\mathcal{P}^*}$ is regular  and hence $\mcr _{\mathcal{P}}$ is regular by \cite[Prop. 1.2.5]{GW2}.
  
  Now suppose that all homogeneous localizations of $R$ at homogeneous primes are regular.  Let $P$ be any prime.  Then $R_{P^*}$ is a regular local ring since it is a localization of $R_{(P^*)}$.  Hence $R_P$ is a regular local ring by \cite[Prop. 1.2.5]{GW2}.  Thus $R$ is a regular ring.
  
  Now suppose $(R, \M)$ is a ${}^*$local ring.  Suppose that $R_{\M}$ is regular.  Let $P$ be any homogeneous prime ideal.  Then, $P \subseteq \M$ implies $R_P$ is regular.  Since $R = R_{(\M)}$ we may deduce that $R$ is regular by the first part of the proof.  The other implication is immediate.
\end{proof}

Now we limit our attention to an affine semigroup ring $\mcr = K[S]$ over a field $K$. We let $\M$ denote the ideal generated by the noninvertible monomials in $\mcr$.

 \begin{notn}
{\rm Notice that for elements  $\alpha, \beta$  of the affine semigroup $S$ the monomial quotient $x^{\alpha}/x^{ \beta} \in \mcr_{\M}$ if and only if  $x^{\alpha}/x^{ \beta} \in \mcr$.  One way to see this is to use the fact that the colon ideal $(x^{\beta}\mcr: x^{\alpha})$ is monomial.  Now let $P$ be a monomial prime of $\mcr$ corresponding to the nonempty face $F$ of $\pos(S)$.  Replacing $S$ by $S_F := S - S \cap F$ and $\mcr$ by $K[S_F] \cong K[S]_{(P)}$, where $K[S]_{(P)}$ denotes the homogeneous localization at $P$, we see that 
$$K[\grp(S)] \cap K[S]_P = K[S]_{(P)}.$$
We will identify $\grp(S \cap F)$  with the  the group of invertible monomials in $\mcr_P$ and refer to $S_F$ as the localization of $S$ at $F$.

Let $S_0$ denote the subgroup of invertible elements in the affine semigroup $S$ and let $\widetilde{S}$ denote the quotient monoid $S/S_0$.  It is well known that  $K[S]$ is regular if and only if $S$ is the direct sum of a free abelian group $\Z^{\ell}$ and a free abelian monoid $\N^k$ (e.g. see \cite[Ex. 6.1.11]{BH}).  
We shall need a slight variant of this result whose proof we omit.}
 \end{notn} 
  
\begin{prop} \label{prop: K[S] regular} The affine semigroup ring $K[S]$ is regular if and only if $\widetilde{S}  \cong \N^k$, where $k = \dim(K[S]_{\M})$.

\end{prop}

Suppose the affine semigroup ring is regular.  Notice that  if  the images of the elements $\gamma_1, \ldots , \gamma_k \in S  \setminus S_0$  form a free basis for $\widetilde{S}$, then the monomials $x^{\gamma_1} , \ldots , x^{\gamma_k}$ form a regular system of parameters for $K[S]_{\M}$.  The local version of the above proposition is the following.
  
  \begin{prop}  Let $P$ be a prime ideal in the affine semigroup ring $\mcr = K[S]$ over a field $K$ corresponding to the face $F$ of $\pos(S)$.    Then, $\mcr_P$ is regular if and only if the quotient  $\widetilde{S_F}$  is free. 
 \end{prop}

  \begin{proof} Notice that the group of units in $S_F$ is $\grp(S \cap F)$.  By Proposition \ref{prop regular graded ring}, $\mcr_P$ is regular if and only if the homogeneous localization $\mcr_{(P)} \cong K[S_F]$ is regular.   The result now follows immediately from Proposition \ref{prop: K[S] regular}. 
  \end{proof}
  
   We now turn our attention  to an alternate characterization of the regularity condition in the spirit of a result in \cite{Vi}.  One advantage of the alternate characterization is that it can be easily checked using the program NORMALIZ \cite{Norm}.  We first prove some auxiliary results.

  \begin{lemma}\label{lemma 1} Let $F = \pos(S) \cap H$ be a face of the positive cone of the affine semigroup $S$ and $\gamma_1 , \ldots , \gamma_k \in S$.    
Suppose that  $\widetilde{S_F}$  is a free abelian monoid and the images of $\gamma_1 , \ldots , \gamma_k $ form a basis.    
    Let $P$ denote the monomial prime of $K[S]$ corresponding to $F$.
  The following assertions hold.
  \begin{enumerate}

  \item  $F$ is contained in precisely $k$ facets $F_i = \pos(S) \cap H_i$ of $\pos(S)$, and hence $F = F_1 \cap \cdots \cap F_k$;
  \item $\sigma_i(\gamma_j) = \delta_{i j}$  for all $1 \le i,j \le k$, where $\sigma_i$ is the primitive linear form corresponding to $ H_i$; and
  \item $\grp(S \cap F) = \grp(S) \cap H_1 \cap \cdots \cap  H_k$.
  \end{enumerate}
  \end{lemma}

  \begin{proof}
  
 (1)  Since $x^{\gamma_1} , \ldots , x^{\gamma_k}$ is a regular system of parameters for $\mcr_P$ we know that $x^{\gamma_i}\mcr_P$ is a height one prime of $\mcr_P$.   Thus there exist facets $F_i$ of $\pos(S)$ corresponding to the height one primes $P_i$ of $\mcr$ such that $x^{\gamma_i}\mcr_P = P_i\mcr_P \, (i = 1, \ldots , k)$.  We have $P_1 + \cdots + P_k = P$ since they are both primes contained in $P$ and they are equal after localizing at $P$.  Suppose $G$ is a facet of $\pos(S)$ containing $F$ and let $Q = P_G$.  Then, $Q\mcr_P = x^{\delta}\mcr$ since $\mcr_P$ is a UFD.  Since $x^{\delta} \in P_1 + \cdots + P_k$ we must have $x^{\delta} \in P_i$ for some $i$.  Hence $Q = P_i$. So $F_1, \ldots , F_k$ are precisely the facets of $\pos(S)$ containing $F$ and $P_1, \ldots , P_k$ are precisely the height one primes contained in $P$.   Thus $F = F_1 \cap \cdots \cap F_k$. 
  
  (2)  By construction, $\sigma_i(\gamma_i) > 0$.  Just suppose   $\sigma_i(\gamma_j) > 0$ for some $j \ne i$.  Then, $x_{\gamma_j} \in P_i$, which implies $P_j\mcr_P \subseteq P_i\mcr_j$, which is absurd.  Hence   $\sigma_i(\gamma_j) = 0$ for $ i \ne j$.  Since $\sigma_i$ is primitive, we must have $\sigma_i(\gamma_i) =1$.
  
  (3) Suppose that $\alpha, \beta \in S$ and $\alpha - \beta \in H_1 \cap \cdots \cap H_k$.  There exist nonnegative integers $a_1, \ldots , a_k, b_1, \ldots, b_k $ such that $\alpha - \sum a_i\gamma_i, \beta - \sum b_i\gamma_i \in \grp(S \cap F)$.   Since $\alpha - \beta \in H_1 \cap \cdots \cap H_k$, we must have $a_i = b_i \; (i = 1 , \ldots , k)$ by (3).  Hence $\alpha - \beta \in \grp(S \cap F)$.  Since the opposite containment is clear we have equality of groups.
  \end{proof}
  
  \begin{lemma}\label{lemma 2}  Let $F$ be a face of $\pos(S)$ that is the intersection of $k$ facets $F_1 = \pos(S) \cap H_1, \ldots , F_k = \pos(S) \cap H_k$ of $\pos(S)$ and let $\sigma_i$ be the primitive linear form associated with $H_i \; (i = 1, \ldots , k)$.  Suppose that
  \begin{enumerate}
  \item  there exist $\gamma_1 , \ldots , \gamma_k \in S$ such that $\sigma_i(\gamma_j) = \delta_{i j} \mbox{ for all } 1 \le i, j \le k$; and
  \item $\grp(S \cap F)  = \grp(S) \cap H_1 \cap \cdots \cap H_k$. 
  \end{enumerate}
  Then, $S_F/U(S_F)$ is a free monoid with basis given by the images of $\gamma_1 , \ldots , \gamma_k$ in the quotient monoid.
  \end{lemma}
  
  \begin{proof}  The proof is straightforward. Suppose $\alpha \in S$ and $\sigma_i(\alpha) = a_i (i = 1, \ldots k)$.  Then $\alpha - (\sum a_i \gamma_i) \in \grp(S) \cap H_1 \cap \cdots \cap H_k = \grp(S \cap F)$ implies the image of $\sum a_i \gamma_i$ in the quotient monoid is $\overline{\alpha}$.  Suppose that $a_i , b_i \in \N$ and $\sum a_i \overline{\gamma_i} = \sum b_i \overline{\gamma_i}$.    Then there exists $\mu \in \grp(S \cap F)$ such that $\sum a_i\gamma_i + \mu = \sum b_i \gamma_i$.  By condition (1) we must have $a_i = b_i \; (i = 1, \ldots , k)$. Hence $\widetilde{S_F}$ is free with the asserted basis.
    \end{proof}

Combining the previous two lemmas we immediately obtain the following characterization.

\begin{theorem}\label{thm: R_k}  An affine semigroup ring $\mcr = K[S]$  satisfies condition $\mathrm{R}_{\ell}$ of Serre if and only if for each positive integer $k \le \ell$ and any face $F$ of $\pos(S)$ such that $\mathrm{ht}(P_F) = k$  there exist facets $F_1, F_2, \ldots, F_k$ of $\mathrm{pos}(S)$  such that $F = F_1 \cap F_2 \cap  \cdots \cap F_k$ and the following conditions hold:

\begin{enumerate}
\item there exist $\bgamma_1, \ldots, \bgamma_k \in S$ such that $\sigma_i(\bgamma_j) = \delta_{i j} \mbox{ for all } 1 \le i, j \le k$;  and
\item $\grp(S \cap F_1 \cap \cdots \cap F_k) = \grp(S) \cap H_1 \cap \cdots \cap H_k$.
\end{enumerate}
\end{theorem}
 
We end this section  with an example of an affine semigroup ring that satisfies condition R$_2$ of Serre but doesn't satisfy condition S$_2$. It was inspired by an example suggested to the author by I. Swanson.

\begin{example}  Suppose $K$ is a field and
consider the semigroup $S$ of $\Z^3$ generated by the following vectors: 
$$(1,0,0), (1,3,0), (1,0,3), (1,1,0), (1,2,0), (1,0,1), (1,0,2), (1,2,1), \mbox{and } (1,1,2).$$

Notice that $\grp(S) = \Z^3$ and that $\pos(S) = H_2^+ \cap H_3^+ \cap H_4^+$, where 
$H_2$ and $H_3$ are the indicated coordinate hyperplanes and $H_4$ is defined by the
 primitive linear form $\sigma$, where $\sigma(a,b,c) = 3a -b -c$.   Thus $\pos(S)$ has 3 
 facets $F_2, F_3, F_4$ and 3 codimension 2 faces  $F_{2 3} = F_2 \cap F_3, F_{2 4} = 
 F_2 \cap F_4, F_{3 4} = F_3 \cap F_4$.  One checks that $\grp(S \cap F_2) = \grp(\{ 
 (1,0,0), (1,0,1) \}) = \grp(S) \cap H_2$ and by symmetry, $\grp(S \cap F_3) = \grp(S) \cap 
 H_3$.  We also have $\grp(S \cap F_4) = \grp(\{ (1,3,0), (1,0,3), (1,2,1), (1,1,2)\}) = \grp(\{ 
 (1,0,3), (0,1,-1)\}) = \grp(S) \cap H_4$.  One can also verify that $\grp(S \cap F_2 \cap
  F_3) = \grp(\{ (1,0,0) \}) = \grp(S) \cap H_2 \cap H_3$, $\grp(S \cap F_2 \cap F_4) = \grp(\{ 
  (1,0,3) \}) = \grp(S) \cap H_2 \cap H_4$, and by symmetry $\grp(S \cap F_3 \cap F_4) = 
  \grp(S) \cap H_3 \cap H_4$.   So the group conditions for the affine semigroup ring 
  $K[S]$ to satisfy R$_2$ are satisfied.   For each codimension 2 face $F_{ij}$ we must 
  produce 2 vectors $\bgamma_i, \bgamma_j$ satisfying $\sigma_i(\bgamma_j) =
   \delta_{i j}$, where $\sigma_2, \sigma_3$ are the coordinate functions and 
   $\sigma_4 = \sigma$ is defined above.  The vectors are given below.
$$\begin{array}{ccc}
(i,j) &\bgamma_i & \bgamma_j \\
(2,3) & (1,1,0) & (1,0,1) \\
(2,4) & (1,1,2) & (1,0,2) \\
(3,4) & (1,2,1) & (1,2,0)
\end{array}$$
Thus condition (1) in Theorem \ref{thm: R_k} is satisfied and we may conclude that $K[S]$ is regular in codimension 2.  However, as we shall now see, $K[S]$ is not normal
 so can't possibly satisfy S$_2$.  
 
 Notice that $(1,1,1) = \frac{1}{3}(1,0,0) + \frac{1}{3}(1,3,0) + \frac{1}{3}(1,0,3) \Rightarrow (1,1,1) \in \pos(S)$.   However, $(1,1,1) \notin S$ and hence $K[S]$  is not normal.

 There is another way to see that $K[S]$ is not normal that is more obvious.  Consider the injective homomorphism $\varphi: \Z^3 \rightarrow \Z^3$ defined by $\varphi(\e_1) = (3,0,0), \\
  \varphi(\e_2) = (-1,1,0), \mbox{ and }\varphi(\e_3) = (-1,0,1)$.   The image of $S$ is the semigroup $\tilde{S}$ generated by the vectors 
  $$(3,0,0), (0,3,0), (0,0,3), (2,1,0), (1,2,0), (2, 0,1), (1,0,2), (0,2,1), (0,1,2).$$
  Notice that we have listed all 3-tuples of non-negative integer whose components sum to 
  3 except (1,1,1).  This isomorphism of semigroups induces an isomorphism of $K[S]$ 
  and $K[\tilde{S}] \cong K[x^3, y^3, z^3, x^2y, xy^2, x^2z, xz^2, y^2z, yz^2]$.  The latter 
  has normalization  \\ $K[x^3, y^3, z^3, xyz, x^2y, xy^2, x^2z, xz^2, y^2z, yz^2]$, which is the
 3rd Veronese subring of \\
 $K[x,y,z]$.  Hence $K[S]$ is not normal.
\end{example}

\section{The Rees Algebras of a Special Class of Monomial Ideals}

We now look at the Rees algebras of a special class of integrally closed monomial ideals.

\begin{notn}\label{blam}
{ \rm Let $\blam = (\lam_1, \ldots, \lam_n)$ be a tuple of positive integers and \\
 $J = J(\blam) = (x_1^{\lam_1}, \ldots, x_n^{\lam_n})$, where the ideal in taken inside the polynomial ring \\$K[\xvec{n}] =: R$, and $I = I(\blam) = \bar{J}$.   Thus $I$ is an integrally closed monomial ideal with minimal monomial reduction $J$.  Let $L = \lcm(\lam_1, \ldots , \lam_n)$, $\omega_i = L/\lam_i \; (i = 1, \ldots , n)$, and $\bomega = (\omega_1 , \ldots , \omega_n)$.  Notice that $L = dw$, where $d = \gcd(\lambda_1 , \ldots , \lambda_n)$.}
\end{notn}

We will characterize those monomial ideals $I(\blam)$ whose Rees algebras satisfy $\mathrm{R}_{\ell}$ for some $\ell < n$.

Observe that the Rees algebra $R[It]$ of a monomial ideal
$I$  can always be identified with an affine  semigroup ring over $K$.  Namely, if $I =
(x^{\bbeta_1}, \ldots, x^{\bbeta_r})$ and  $S(I) = \langle (\e_1,0),
\ldots, (\e_n,0),(\bbeta_1,1), \ldots, (\bbeta_r,1) \rangle \subseteq
\N^{n+1}$, then $R[It] \cong K[S(I)]$. 

In case $\mcr := R[It]$ is the Rees algebra of $I = I(\blam)$ the condition that every  height $k$ monomial prime corresponds to an intersection of precisely $k$ facets  and condition (2) of Theorem \ref{thm: R_k} automatically hold as we shall see below.  First we describe the height $k$ monomial primes of $\mcr$. 

 The facets $F_{\sigma}, F_1, \ldots , F_{n+1}$ of
$\pos(S)$ are cut out by the supporting hyperplanes \\ $H_{\sigma}, H_1,
\ldots, H_{n+1}$  where $\sigma(\balpha, a_{n+1})=\bomega \cdot \balpha
- La_{n+1}$ and $ H_1,
\ldots, H_{n+1}$ are the coordinate hyperplanes in $\R^{n+1}$.  Notice that with this notation the generating set for $S(I)$ is
$$ \{ (a_1,\ldots,a_n,d) \in \mathbb{N}^{n+1} \mid
a_1
\omega_1 + \cdots + a_n \omega_n \ge dL \mbox{ for } d \le 1 \}.$$
Notice that $(\e_1,0), \ldots , (\e_n,0), (\lam_1\e_1,1) \in S$ and hence $\pos(S) = \Z^{n+1}$, i.e., $S$ is full-dimensional.

The following description of the height one monomial primes of
$R[It]$ appeared in \cite{Vi}.

\begin{lemma}\label{lem: ht 1 monl prime}  For a monomial ideal $I=I(\blam)$ the height one
monomial primes of $R[It]$ are as follows:
\begin{eqnarray*}
P_i  &= & (x_i) + (x^{\bbeta_j}t \mid \e_i \le_{pr} \bbeta_j) \; for \; (i=1,
\ldots, n);\\
P_{n+1} &=& (x^{\bbeta_1}t, \ldots, x^{\bbeta_r}t); and  \\
 P_{\sigma} &= & (x_1, \ldots, x_n) + (x^{\bbeta_j}t \mid
\sigma(\bbeta_j,1) > 0).
\end{eqnarray*}
\end{lemma}

We wish to describe the height $k$ monomials ideals.  Towards this end we show that every codimension $k$ face of $\pos(S)$ is the intersection of precisely $k$ facets of $\pos(S)$ for each $k$ such that $1 \le k \le n$.   Notice that 
\begin{eqnarray*}
\pos(S(I)) &=& \pos(S(J)) \\
&=&  \pos((\e_1,0), \ldots, (\e_n,0), (\lam_1\e_1,1), \ldots , (\lam_n\e_n,1)).
\end{eqnarray*}

  Alternate proofs that the following are the height $k$ monomial primes of $R[It]$ can be found in \cite{Co}.

In the next few paragraphs, given integers $1 \le i_1 < i_2 < \cdots < i_k \le n$ let $1 \le j_1< \cdots < j_{n-k} \le n$ be such that $\{i_1, \ldots , i_k, j_1, \ldots , j_{n-k} \} = \{ 1 , \ldots , n \}$.

\begin{lemma}\label{lem: face 1} For $k < n$ and integers $1 \le i_1 < i_2 < \cdots < i_k \le n$, let \\ $F= F_{i_1} \cap \cdots \cap F_{i_k}$.  Then,
\begin{enumerate}
\item $ F = \pos((\e_{j_1},0), \ldots , (\e_{j_{n-k}},0), (\lam_{j_1}\e_{j_1},1), \ldots , (\lam_{j_{n-k}}\e_{j_{n-k}},1));$
\item  $\codim(F_{i_1} \cap \cdots \cap F_{i_k}) = k$ and 
$P_F = (x_{i_1}, \ldots , x_{i_k}) + (x^{\bbeta}t \mid  (\bbeta ,1) \in S \setminus F);$ and 
\item $\widetilde{S_F}$ is free with basis given by the images of $(\e_{i_1},0), \ldots , (\e_{i_k}, 0)$ and hence $R[It]$ localized at $P_F$ is regular.
\end{enumerate} 
\end{lemma}

\begin{proof} Let $\sigma_i(\balpha, a_{n+1}) = a_i$ for $1 \le i \le n+1$. 
The first assertion is a consequence of the fact that $\sigma_i$ of a sum of vectors in $\pos(S)$ is zero if and only if  $\sigma_i$ of each summand is zero.  The codimension statement follows immediately.  The description of the corresponding monomial prime comes from looking at the generators of the prime ideal $S \setminus F$ of $S$. 

To see that the images of  $(\e_{i_1},0), \ldots , (\e_{i_k}, 0)$  generate the quotient monoid it suffices to consider generators of $S$ of the form $(\beta,1)$ that aren't in $S \cap F$.  For such, $b_{i_s} > 0 $ for some $s$.  
 
Then, 
\begin{eqnarray*}
 (\beta,1)  &\cong &(( \beta ,0) - (\lam_{j_1} \e_{j_1}, 0)) +  (\lam_{j_1} \e_{j_1}, 1) \pmod{\grp(S \cap F)}  \\
 &\cong& (b_{i_1}\e_{i_1},0)  + \cdots + (b_{i_k}\e_{i_k},0) \hspace{.35in} \pmod{ \grp(S \cap F)}
 \end{eqnarray*}
Thus the images of $(\e_{i_1},0), \ldots , (\e_{i_k}, 0)$ form a free basis for the quotient  $\widetilde{S_F}$, since uniqueness of representation is clear.
\end{proof}

Next we consider the case where one of facets in the intersection is $F_{n+1}$. The proof is due to the same observations that appeared in  the preceding proof, so we omit it.

\begin{lemma} \label{lem: face 2} For  integers $1 \le i_1 < i_2 < \cdots < i_{k} \le n$ let $F=F_{i_1} \cap \cdots \cap  F_{i_{k}} \cap F_{n+1}$.  The following hold.   
\begin{enumerate}
\item $F=  \pos((\e_{j_1}, 0) , \ldots, (\e_{j_k}, 0)).$ 
\item $\codim(F) = k+1$ and the corresponding monomial prime is \\
$P_F = (x_{i_1}, \ldots , x_{i_k}) + (x^{\bbeta}t \mid (\bbeta, 1) \in S).$
\item If $k < n$ then, $\widetilde{S_F}$ is free with basis given by the images of the vectors \\ $(\e_{i_1},0), \ldots , (\e_{i_k},0), (\e_{j_1},1)$. In case $k < n$ the Rees algebra $R[It]$ localized at $P_F$ is regular.
\end{enumerate}
\end{lemma}

Notice that $F_1 \cap \cdots \cap F_n = \{(0, \ldots , 0)\} = F_{n+1} \cap F_{\sigma}$ and the corresponding monomial prime is $\M = (x_1 , \ldots , x_n) + (x^{\bbeta}t \mid (\bbeta , 1) \in S)$, the maximal monomial prime of $K[S]$.   In  this case, the codimension drops more than the expected amount. We can still realize the apex of the cone as the intersection of $n+1$ facets, namely $ \{(0, \ldots , 0)\} = F_1 \cap \cdots \cap F_{n+1}$.  Notice $S_0 = \{ 0 \}$ and $\widetilde{S}$ is not free provided that $n > 1$.

Now we involve the facet $F_{\sigma}$.  By the above lemmas these are the only faces we need to worry about when characterizing which Rees algebras $R[It]$ are regular in codimension $k \le n$.  The proof of the next lemma is a consequence of the same observations and is omitted. 

\begin{lemma} \label{lem: face 3} For integers $1 \le i_1 < i_2 < \cdots < i_k  \le  n$ let  $F = F_{i_1} \cap \cdots \cap F_{i_k} \cap F_{\sigma}$.  The following hold.
\begin{enumerate}
\item $F =  \pos((\lam_{j_1}\e_{j_1},1), \ldots , (\lam_{j_{n-k}}\e_{j_{n-k}},1)).$ 
\item $\codim(F_{i_1} \cap \cdots \cap F_{i_k} \cap F_{\sigma} ) = k+1$ and the corresponding monomial prime is
$P_F = (x_{i_1}, \ldots , x_{i_k}) + (x^{\bbeta}t \mid (\bbeta, 1) \in S \setminus F).$
\end{enumerate}
\end{lemma}

We now show that
condition (2) of Theorem \ref{thm: R_k} is always satisfied by the positive dimensional faces of $\pos(S)$ that are contained in $F_{\sigma}$.  The condition is a priori true for positive dimensional faces not contained in $F_{\sigma}$ by Lemmas  \ref{lem: face 1} and \ref{lem: face 2}.

\begin{lemma}\label{lem: ker Z^n to Z} Let  $\bomega = (\omvec{n})$ be a tuple of positive integers with $n \ge 2$  and let 
$\phi: \Z ^n \rightarrow \Z $ be defined by $\phi(\balpha) = \bomega \cdot \balpha$.  For each $1 \le i < j \le n$ let $r_{ij} = \gcd(\omega_i, \omega_j)$.  Then, the kernel of $\phi$ is generated by the  tuples $\bmu_{ij} = \frac{\omega_j}{r_{ij}}\e_i - \frac{\omega_i}{r_{ij}}\e_j \; \; (1 \le i < j \le n)$, where $\e_1, \ldots , \e_n$ are the standard basis vectors for $\Z^n$. Furthermore, for integers $1 \le i_1 < i_2 < \cdots < i_k  \le  n$ we have $\ker(\phi) \cap H_{i_1} \cap \cdots \cap H_{i_k}$ is generated by vectors $\bmu_{ij}$ in $H_{i_1} \cap \cdots \cap H_{i_k}$.
\end{lemma}
\begin{proof}  This lemma is used in Gr$\mathrm {\ddot{o}}$bner basis theory but for the reader's convenience we will supply a proof.  

We proceed by induction on $n$ the case $n=2$ being straightforward.  Suppose $n > 2$ and assertion holds for $n - 1$.  Assume $\bbeta = (b_1, \ldots, b_n) \in \ker(\phi)$. Let $g = \gcd(\omega_1, \ldots, \omega_n)$.  Then, \begin{equation}
\label{ 1}
b_n\frac{\omega_n}{g}= - (b_1\frac{\omega_1}{g} + \cdots + b_{n-1}\frac{\omega_{n-1}}{g}),  
\end{equation} 
which implies
\begin{equation}
\label{2 }
b_n =  \sum_{i=1}^{n-1} c_i\frac{\omega_i}{g} =  \sum_{i=1}^{n-1} c_is_i\frac{\omega_i}{r_{in}},
\end{equation}
where $r_{in} = s_ig \, (i=1, \ldots, n-1)$.  Then,
\begin{equation}
\label{ 3}
\bbeta + \sum_{i=1}^{n-1} c_is_i\bmu_{in} = (b_1^{\prime}, \ldots, b_{n-1}^{\prime},0) \in \ker(\phi),
\end{equation}
and the assertion then follows from the induction hypothesis.  Notice that if $\bbeta \in \ker(\phi) \cap H_{i_1} \cap \cdots \cap H_{i_k}$ then each step only involves vectors in $H_{i_1} \cap \cdots \cap H_{i_k}$.
\end{proof}

This enables us to prove that in our setting the group property is automatic.  The following two results appear in the unpublished thesis of H. Coughlin \cite{Co}.

\begin{lemma}  \label{lem: group 2} For integers $1 \le i_1 < i_2 < \cdots < i_k  \le  n$ we always have
$$\grp(S \cap F_{i_1} \cap \cdots \cap F_{i_k} \cap F_{\sigma} ) = \grp(S) \cap H_{i_1} \cap \cdots \cap H_{i_k} \cap H_{\sigma} .$$
\end{lemma}

\begin{proof}  Recall that  $\grp(S) = \Z^{n+1}$.  The containment $\grp(S \cap F_{i_1} \cap \cdots \cap F_{i_k} \cap F_{\sigma} ) \subseteq \grp(S) \cap H_{i_1} \cap \cdots \cap H_{i_k} \cap H_{\sigma} = \Z^{n+1}\cap H_{i_1} \cap \cdots \cap H_{i_k} \cap H_{\sigma} $ is clear.  

If $k = n$ then $\grp(S) \cap H_{i_1} \cap \cdots \cap H_{i_k} \cap H_{\sigma} = \Z^{n+1} \cap H_1 \cap \cdots \cap H_n \cap H_{\sigma} =  \{ \mathbf{0} \}$ and the assertion follows.  Now assume $k < n$.  Suppose $(\bbeta,d) \in \Z^{n+1}\cap H_{i_1} \cap \cdots \cap H_{i_k} \cap H_{\sigma} $.  Then, $\bbeta - d\lam_n\e_n \in \ker(\phi)$, where $\phi$ is defined as in the preceding lemma.  By that lemma and the fact that $(\lam_{j_1}\e_{j_1},1) \in S \cap F_{i_1} \cap \cdots \cap F_{i_k} \cap F_{\sigma} $, it suffices to show that $(\bmu_{ij},0) \in \grp(S \cap H_{i_1} \cap \cdots \cap H_{i_k} \cap H_{\sigma})$ for $1 \le i < j \le n$ and $\bmu_{i j} \in H_{i_1} \cap \cdots \cap H_{i_k}$, where we are viewing $H_{i_j}$ as a coordinate hyperplane in either $\R^n$ or $\R^{n+1}$.  

Let $1 \le i < j \le n$, set $r = \gcd(\omega_i,\omega_j)$ and choose $d \in \Z$ such that $0 \le d\lam_j - \frac{\omega_i}{r} < \lam_j$.  Multiplying by $\omega_j$ and dividing by $\omega_i$ we get  $0 \le d\lam_i - \frac{\omega_j}{r}  < \lam_i$ which implies $0 <   \frac{\omega_j}{r} - (d-1)\lam_i  \le \lam_i$.   Then,
\begin{eqnarray*}
(\bmu_{ij},0) &=& (\frac{\omega_j}{r}\e_i +(d\lam_j - \frac{\omega_i}{r})\e_j,d) - d(\lam_j\e_j,1)  \\
& = & ((\frac{\omega_j}{r} - (d-1)\lam_i)\e_i + (d\lam_j - \frac{\omega_i}{r})\e_j,1) + (d-1)(\lam_i\e_i,1)-d(\lam_j\e_j,1)  \\
& \in & \grp(S \cap F_{i_1} \cap \cdots \cap F_{i_k} \cap F_{\sigma}), 
\end{eqnarray*}
since each of the 3 tuples involved is a generator of $S$ that is in $F_{i_1} \cap \cdots \cap F_{i_k} \cap F_{\sigma}$.
\end{proof}

Combining Theorem \ref{thm: R_k} with the preceding lemma we obtain the following result.

\begin{theorem} \label{thm: R_k for Rees} For a positive integer $\ell < n$ the Rees 
algebra of $I = I(\blam)$ over a field satisfies condition $\mathrm{R}_{\ell + 1}$ of 
Serre if and only if for all sequences of positive integers $1 \le i_1 < \cdots < i_{\ell} \le n$
there exist $\bgamma_i = (\bbeta_i , 1) \in \N^{n+1} \; (i = 1 , \ldots , \ell +1)$ such that $\sigma_i(\bgamma_j) = 
\delta_{i j} \, (1 \le i,j \le \ell + 1)$, where $\sigma_1 , \ldots , \sigma_{\ell +1}$ are the primitive linear forms associated with the hyperplanes $H_{i_1} , \ldots , H_{i_{\ell}}, H_{\sigma}$.  
\end{theorem}

\begin{proof}  First let us determine when $R_{\ell + 1}$ holds. We need only consider faces contained in $F_{\sigma}$. By Lemma \ref{lem: group 2} and
 Theorem \ref{thm: R_k} we need only show that condition (1) of Theorem \ref{thm: R_k}
  holds for such faces.  Let $1 \le k \le \ell +1$ and let $Q$ be a height $k$ monomial prime corresponding to
  a face contained in the facet $F_{\sigma}$.  Then $Q$ is contained in a height $\ell +1$
   monomial prime corresponding to a face contained in the facet $F_{\sigma}$ by Lemma
   \ref{lem: face 3}.  Thus it suffices to establish condition (1) of Theorem \ref{thm: R_k} for
    height $\ell+1$ monomial primes whose faces are contained in the facet $F_{\sigma}$. 
Hence it is necessary and sufficient that there exist  vectors $\bgamma_i  \in S(I) \; (i = 1 , \ldots , \ell +1)$ such that
    $\sigma_i(\bgamma_j) = \delta_{i j}$ where $\sigma_1 , \ldots , \sigma_{\ell +1}$ are the
     primitive linear forms associated with the hyperplanes $H_{i_1} , \ldots , H_{i_{\ell}},
      H_{\sigma}$.   
      
      Recall that the generators of $S(I)$ have $(n+1)^{\mathrm{st}}$ component 0 or 1. Write each $\bgamma_j$ as a sum of generators of $S(I)$. First suppose that $1 \le j \le \ell $. The condition
      that $\sigma_i(\bgamma_j) = \delta_{i j} \; (i = 1, \ldots , \ell + 1)$ forces  some summand $(\bbeta_j, 1)$ of 
      $\bgamma_j$ to satisfy $\sigma_i(\bbeta_j,1) = \delta_{i j}$ for all $1 \le i \le \ell +1$.  Replacing $\bgamma_j$ by this summand we may assume that  $\bgamma_j = (\bbeta_j,1).$ 
      Now 
      consider the summands involved in the expression for $\bgamma_{\ell + 1}$. Consider 
      first the possibility that all summands have $(n+1)^{\mathrm{st}}$ component 0.  Then,
       each summand has positive $\sigma$-value so there is only  one summand 
       $\bgamma_{\ell +1} = (\e_j ,0)$, where $j \in \{j_1, \ldots , j_{n_\ell}  \}$ and 
       $\omega_j = 1$.  In this case we also have $( (L+1)\e_{j} ,1)$ satisfies the requirements and we may replace $\bgamma_{\ell+1}$ by $((L+1)\e_{j} ,1)$.  The remaining possibility is that some summand has $(n+1)^{\mathrm{st}}$ component 1 and again we may replace $\bgamma_{\ell+1}$ by this summand.  In any case, we may assume $\bgamma_{\ell + 1}$ has $(n+1)^{\mathrm{st}}$ component 1.     Conversely, if vectors $(\bbeta_j, 1) \in \N^{n+1}$ satisfying $\sigma_i(\bbeta_j,1) = \delta_{i j}$  for all $1 \le i, j \le \ell +1$ exist, they are automatically in $S(I)$ since $I$ is integrally closed. 
      \end{proof}

We now state the result entirely in terms of the integers $L, \omega_1, \ldots , \omega_n$ determined by the vector  $\blam = (\lam_1 , \ldots , \lam_n)$.

\begin{cor}  \label{cor: R_k for Rees} For a positive integer $\ell < n$ the Rees 
algebra of $ I(\blam)$ over a field satisfies condition $\mathrm{R}_{\ell + 1}$ of 
Serre if and only if for all sequences of positive integers $1 \le i_1 < \cdots < i_{\ell} \le n$
 and $1 \le j_1 < \cdots < j_{n-\ell} \le n$ such that $\{1, \ldots , n \} = \{i_1 , \ldots , i_{\ell} \} \cup \{ j_1 , \ldots , j_{n-\ell}
  \} $, we  have $$L - \omega_{i_1}, \ldots , L - \omega_{i_{\ell}}, L + 1 \in \langle \omega_{j_1}, \ldots , \omega_{j_{n-\ell}} \rangle.$$
\end{cor}

\begin{proof}  By Theorem \ref{thm: R_k for Rees} it suffices to show that for a sequence of positive integers $1 \le i_1 < \cdots < i_{\ell} \le n$ 
 and $1 \le j_1 < \cdots < j_{n-\ell} \le n$ such that $ \{1, \ldots , n \} = \{i_1 , \ldots , i_{\ell} \} \cup \{ j_1 , \ldots , j_{n-\ell}  \} $, 
  vectors $\bgamma_j = (\bbeta_j, 1) \; (j = 1, \ldots, \ell +1) \in \N^{n+1}$ exist such that $\sigma_i(\bgamma_j) = \delta_{i j} \;  \mbox{ for all }1 \le  i, j \le \ell + 1$, where the $\sigma_i$ are as above, if and only if 
$L - \omega_{i_1}, \ldots , L - \omega_{i_{\ell}}, L + 1 \in \langle \omega_{j_1}, \ldots , \omega_{j_{n-\ell }} \rangle.$
         
         Suppose first that the vectors $\bgamma_i = (\bbeta_i, 1) \in \N^{n+1}$ satisfying the
          necessary conditions exist.  By our requirements, $\bgamma_i= (\e_i + a_{j_1}\e_{j_1} + \cdots + a_{j_{n - \ell}}\e_{j_{n - \ell}},1)$ for all $i = 1, \ldots , \ell$, where the coefficients $a_{j_1} , \ldots , a_{j_{n-\ell}}$ are nonnegative.   The existence of 
          such vectors $\bgamma_i$ is equivalent to the conditions $L - \omega_{i_1}, \ldots , L - \omega_{i_{\ell}} \in \langle \omega_{j_1}, \ldots , \omega_{j_{n-\ell}} \rangle.$  We
           must also have $\bgamma_{\ell + 1} = (a_{j_1}\e_{j_1} + \cdots + a_{j_{n - \ell}}\e_{j_{n - \ell}} ,1)$, where the coefficients $a_{j_1} , \ldots , a_{j_{n-\ell}}$ are
            nonnegative and  $ a_{j_1}\omega_{j_1} + \cdots + a_{j_{n - \ell}}\omega_{j_{n - \ell}}  - L = 1$.  The existence of such a $\bgamma_{\ell + 1}$ is equivalent to   $L + 1 \in 
            \langle \omega_{j_1}, \ldots , \omega_{j_{n-\ell }} \rangle.$
\end{proof}

Applying this theorem for values of $\ell$ close to $n$ gives simple descriptions of when the Rees algebra of $I(\blam)$ is regular in codimension $\ell$.

\begin{cor}  The Rees algebra of $ I(\blam) \subset K[\xvec{n}]$ over a  field $K$ is regular in codimension $n$ if and only if $\blam = \lam(1,1, \ldots , 1)$  and hence, $I(\blam) = \M^{\lam}$.
\end{cor}

\begin{proof}  The Rees algebra of $I$ satisfies $\mathrm{R}_{n}$ if and only if for every  sequence  \\$1 < \cdots < i-1 < i+1 < \cdots < n$ of length $n-1$, we have
$$L - \omega_1, \ldots , L - \omega_{i-1} , L - \omega_{i+1} , \ldots , L-\omega_n , L+1 \in \langle \omega_i \rangle.$$  In particular, $L+1 = a\omega_i$ for some $a \ge 0$, which implies $1 = \omega_i(a - \lam_i)$.  Thus each $\omega_i = 1$.  Conversely, if all the $\omega_i = 1$ then the necessary conditions are satisfied.
\end{proof}

\begin{cor} \label{cor: 2} Suppose that $n \ge 3$. The Rees algebra of $ I(\blam) \subset K[\xvec{n}]$ over a field $K$ is regular in codimension $n-1$ if and only if the integers $\omega_i$ are pairwise relatively prime.
\end{cor}

\begin{proof}  The sequences of length $n-2$ arise from omitting two integers $1 \le i < j \le n$.  For each pair $1 \le i < j \le n$ we must have $$L - \omega_k \in \langle \omega_i , \omega_j \rangle \mbox{ for all } k \ne i,j \mbox{ and } L+1 \in \langle \omega_i , \omega_j \rangle.$$  Write $L+1 = a\omega_i + b\omega_j$ for $a, b \ge 0$ and read modulo 
$\omega_i$ to obtain the congruence $b\omega_j \equiv 1 \pmod{\omega_i}$.  Hence $\omega_i$ and
 $\omega_j$ are relatively prime.  This holds for all pairs $1 \le i<j \le n$. Conversely, if the 
 integers $\omega_i$ are pairwise relatively prime then every integer at least $(\omega_i - 1)(\omega_j -1)$ is in $\langle \omega_i,  \omega_j \rangle$. So $L + 1 \mbox { and } L - \omega_k = \omega_k(g \prod_{s \ne k} \omega_s  -  1) \in \langle \omega_i,  \omega_j \rangle$, where $g = \gcd(\lam_1, \ldots , \lam_n)$.
\end{proof}

If $n = 2$ the Rees algebra $R[I(\blam)t]$ is normal and hence regular in codimension 
$n - 1 = 1$ without any additional assumptions. Corollary \ref{cor: 2} says that if $n=3$ 
the Rees algebra of $I(\blam)$ is regular in codimension 2 if and only if the $\omega_i$ 
are pairwise relatively prime.  H. Coughlin \cite{Co} proved this special case and also that this condition is sufficient for $R[I(\blam)t]$ to be normal.   This result combined with an earlier result of Reid, Roberts, and Vitulli \cite{RRV}  has an interesting consequence, which we now present. 

\begin{notn}
{\rm For an $\N$-graded ring $A$ and a positive integer $t$ we let $A_{\ge t}$ denote the homogeneous ideal $A_{\ge t} = \oplus_{s \ge t} A_s$.  Further assume that $A$ is generated as an $A_0$-algebra by homogeneous elements $\xvec{n}$ of positive degrees $\omega_1 , \ldots , \omega_n$, respectively.

Notice that for a tuple of positive integers $\blam$, with $L = \lcm(\lambda_1, \ldots, \lambda_n )$ and $\omega_i = L/\lambda_i  \,  (i = 1 , \ldots , n)$ as in (\ref{blam}), if we define a new grading on  $R = K[\xvec{n}]$ by declaring $\deg(x_i) = \omega_i \,  (i = 1 , \ldots , n)$, then $I(\blam) = R_{\ge L} = R_{\ge dw}$, where $d = \gcd(\lam_1, \ldots , \lam _n)$.  Ideals of the form $A_{\ge t}$ have been studied by E. Hyry and K. E.  Smith  in connection with Kawamata's Conjecture which speculates that every adjoint ample line bundle on a smooth variety admits a nonzero section (e.g. see \cite{HS1} and \cite{HS2}).  They also arise as test ideals in tight closure theory as illustrated in \cite[Remark 3]{F}}.
\end{notn}
 
 \begin{prop}  Let $R= K[x,y,z]$ be a polynomial ring over a field $K$ and let $a,b,c$ be pairwise relatively prime positive integers.  Set $S = K[x^a, y^b, z^c]$ and $L = abc$.  Then, the ideal $\mathfrak{a} = S_{\ge L}$ is normal, that is,  $\mathfrak{a}^t = S_{\ge tL}$ for all $t \ge 1$.
 \end{prop}
 
 \begin{proof}  Observe that the homogeneous ideal $ \mathfrak{a} = S_{\ge L}$ is 
 integrally closed and that the integral closure of $ \mathfrak{a}^t$ is   $ S_{\ge t L}$  for $t \ge 1$
 (e.g. see the discussion in \cite{RRV}). By Proposition 3.7 of  \cite{RRV} to prove that $\mathfrak{a}$ is normal,  it suffices to show that $\mathfrak{a}^2 = S_{\ge 2L}$.  For this we proceed as in the proof of  Theorem III.2.2 of \cite{Co}.  Suppose that $x^{ua}y^{bv}z^{cw} \in S_{\ge 2L}$ is a minimal monomial generator.  In particular,  $ua + vb + wc \ge 2L$.  We must exhibit a decomposition $(u,v,w) = (u_1,v_1,w_1)+(u_2,v_2,w_2)$ with $u_ia + v_ib + w_ic \ge L \, (i=1,2)$.  If $u \ge L/a, v  \ge  L/b$ or $w  \ge L/c$ it is clear that we can do this.  For example, if $u \ge L/a$ write $(u,v,w) = (L/a,0,0) + (u - L/a,v,w)$. Thus it suffices to assume that $u < L/a, v <  L/b$ and $w < L/c$.  Notice that this forces the sum of any two of $ua, vb,$ and $wc$ to be strictly greater than $L$ and each summand to be positive.  We consider three cases.
 
First suppose that either $a, b,$ or $c$ is 1. For example, suppose $a=1$. 
 Say $u + vb = L + u_2$.  Then $0 < u_2 < u$ and $$(u,v,w) = (u-u_2,v,0) + (u_2, 0, w)$$ is the desired decomposition.
Thus we may assume that $a,b,c > 1$.  
 
Now suppose that $u <  L/2a , v <  L/2b  $ or $w < L/2c $.  Without loss of generality we may assume that $w  <  L/2c $.  Then, $L - wc > L/2 \ge (a-1)(b-1)$ so there exist $u_1, v_1 \in \N$ such that $u_1a + v_1b = L - wc$.  Since $vb < L$, we have $ua + wc > L$.  Now $u_1a \le L - wc < ua$ implies $u_1 < u$.  Similarly, $v_1 < v$.  Therefore, $$(u,v,w) = (u_1,v_1,w) + (u-u_1, v-v_1,0)$$ is the desired decomposition.

Finally, assume that $u \ge L/2a , v \ge  L/2b  $ and $w  \ge L/2c $.  Set $w_1 = \lceil L/2c  \rceil$.  Then, 
$L - w_1c \ge L - (ab+1)c/2 = (c/2)(ab-1) > (a-1)(b-1)$ and we may write $u_1a + v_1b  = L - w_1c$ for some $u_1, v_1 \in \N$.  Notice that $u_1a \le L - w_1c \le L/2$ and hence $u_1 \le  L/2a  \le  u$.  Similarly, $v_1 \le  v$.  Thus $$(u,v,w) = (u_1,v_1,w_1) + (u-u_1,v-v_1,w-w_1)$$ is the desired decomposition.
 \end{proof}
 
We now present an example of a Rees algebra $\mcr = R[It]$ of a monomial ideal that satisfies $\mathrm{R}_2$ but is not normal.   In order to find an example we must work over a polynomial ring in at least 4 indeterminates by the above remarks.  The following  example  is due to H. Coughlin \cite{Co}.   The example was first explored using the program NORMALIZ \cite{Norm} of Bruns and Koch.

\begin{example}\label{ex: weirdex} Let $\blam=(1443,37,21,91)$. Define $I=I(\blam)$, $S=S(I)$, and $\mcr$ as above. We claim $\mcr$ is not normal but satisfies the equivalent conditions for R$_2$. Hence,  $\mcr$ does not satisfy S$_2$.
\end{example}

We first show that $\mcr$ is not normal
Note that $L=10101$.
  The vector $\balpha=(2,36,1,89)$ satisfies 
$\bomega  \cdot  \balpha=2L$, so that $x^{\balpha}  \in \overline{I^2}$.   Direct computation shows that   $\balpha$ is not the sum of two vectors in $S$ and hence $x^{\balpha} \notin I^2$.   Thus $\mcr$ is not normal.

We show that $\mathrm{R}_2$ holds.  As in the proof of Theorem \ref{thm: R_k for Rees} we need only deal with the height two monomial primes corresponding to the 
faces $G_1, G_2, G_3, G_4$, where $G_i = \pos(S) \cap H_i \cap H_{\sigma}$.  As in the proof of
Theorem \ref{thm: R_k for Rees}, for each $G_i$ we must define  a 
pair of elements $\bgamma_i$ and $\bgamma_6$ in $S$ such that $\sigma_a(\bgamma_b) = \delta_{a b}$ for $a, b \in \{ i, 6 \}$. 

The following vectors satisfy $\sigma_a(\bgamma_b)=\delta_{ab}$ for $a, b \in \{ i, 6 \}$:
$$\begin{array}{ccc}
i&\bgamma_i & \bgamma_6\\ 
1&(1,28,8,12,1) &(0,8,1,67,1) \\

2&(35,1,1,82,1)&(16,0,0,91,1)\\
3&(35,1,1,82,1) & (16,0,0,91,1)\\
4&(220,1,17,1,1) & (275,0,17,0,1)\\
\end{array}$$

$\mcr$ satisfies R$_2$ by Theorem~\ref{thm: R_k for Rees}.  Thus $\mcr$ does not satisfy S$_2$.



\begin{thebibliography}{40}

\bibitem{BG}
W. Bruns and J. Gubeladze,
{Rectangular Simplicial Complexes},
in:  D. Eisenbud (Ed.), { Commutative Algebra, Algebraic Geometry, and Computational Methods},
Springer-Verlag,  Singapore, 1999, 201--214.

\bibitem{BGT}
W. Bruns and J. Gubeladze and N.V. Trung,
{Normal polytopes, triangulations, and Koszul Algebras},
{J. Reine Angew. Math.} {\bf 484} (1997), 123--160.

\bibitem{BH}
W. Bruns and J. Herzog,
{Cohen-Macaulay Rings},
Cambridge University Press, 1993.

\bibitem{Norm}
W. Bruns and R. Koch,
{NORMALIZ: a Program to Compute Normalizations of Semigroups. Available by anonymous ftp from \textbf{ftp://ftp.mathematik.UniOsnabrueck.DE/pub/osm/kommalg/software/
}}.

\bibitem{Co} H. Coughlin, Classes of normal monomial ideals, University of Oregon Ph.D. Dissertation, August 2004.


\bibitem{Ew}
G. Ewald,
{Combinatorial Convexity and Algebraic Geometry}, Graduate Texts in Mathematics {\bf 168}, Springer, New York, 1996.

\bibitem{F} 
S. Faridi,
{Normal ideals of graded rings},
{ Comm. Algebra} {\bf 28} (2000), no. 4, 1971--1977

\bibitem{Gil}
R. Gilmer
{Commutative Semigroup Rings}, Chicago Lectures in Mathematics, The University of Chicago Press, Chicago and London, 1984.

\bibitem{GSW76}
S. Goto and N. Suzuki and K. Watanabe,
{On affine semigroup rings},
{Japan. J. Math.} {\bf 2} (1976), No. 1, 1--12. 

\bibitem{GW1}
S. Goto and K. Watanabe,
{On Graded Rings, I},
{J. Math. Soc. Japan} {\bf 30} (1978), 97--154.

\bibitem{GW2}
S. Goto and K. Watanabe,
{On Graded Rings, II ($Z^n$-graded rings)}, 
{Tokyo J. Math.} {\bf 1} (1978), No. 2, 237--261.

\bibitem{Grun}
B. Gr\"{u}nbaum, Convex Polytopes, Second Edition, {Graduate Texts in Mathematics}  {\bf 221}, Springer, New York, 2003.

\bibitem{HT}
L.T. Hoa and N.T. Trung,
{Affine semigroups and Cohen-Macaulay rings generated by monomials},
{Trans. of Amer. Math. Soc. }   {\bf 298} (1986), No. 1, 145--167. 

\bibitem{HT2}
L.T. Hoa and N.T. Trung,
{Corrigendum to ``Affine semigroups and Cohen-Macaulay rings generated by monomials"},
{Trans. of Amer. Math. Soc. }  {\bf 305} (1988) No. 2, 857. 


\bibitem{Hoch}
M. Hochster,
{Rings of invariants of tori, Cohen--Macaulay rings generated by monomials and polytopes},
{ Ann. of Math.} {\bf 96} (1972), 318--337.

\bibitem{HS1}
E. Hyry, and K.E. Smith,
{On a non-vanishing conjecture of Kawamata and the core of an ideal},
{Amer. J. Math.} {\bf 125} (2003), no. 6, 1349--1410. 

\bibitem{HS2}
E. Hyry and K.E. Smith,
{Core versus graded core, and global sections of line bundles},
{ Trans. of Amer. Math. Soc.} {\bf 356} (2004), no. 8, 3143--3166.
 
\bibitem{Ish}
M.N. Ishida,
{The local cohomology groups of an affine semigroup ring}, in: H. Hijikata et. al. (Eds.), 
{ Algebraic geometry and commutative algebra in honor of Masayoshi Nagoya}, {\bf Vol. 1},  
 Kinokuniya, Tokyo (1988), 141--153.


\bibitem{RRV}
L. Reid and L.G. Roberts and M.A. Vitulli,
{Some results on normal homogeneous ideals},
{ Comm. in Alg}, 
\textbf{31} (2003), No. 9, 4485--4506.



\bibitem{Vi}
M.A. Vitulli,
{On normal monomial ideals}, in:  C.G. Melles, J-P Brasselet, G. Kennedy, K. Lauter, L. McEwan, (Eds.), Topics in Algebraic and Noncommutative Geometry,
{Contemporary Mathematics} {\bf 324} (2003), 205--218.

\bibitem{V}
M.A. Vitulli,
{Weak normalization and weak subintegral closure}, in: C.G. Melles and R.I. Michler (Eds.), {Singularities in Algebraic and Analytic Geometry},
 {Contemporary Mathematics} {\bf 266} (2000), 11--21.


\bibitem{Zieg} G. M. Ziegler, { Lectures on Polytopes}, Graduate Texts in Mathematics {\bf 152}, Springer, New York 1998.


\end{thebibliography}
\end{document}